\documentclass[onetabnum,letterpaper]{siamltex}
\usepackage{amsfonts,amsmath}
\usepackage{subfig}
\usepackage{MnSymbol}
\usepackage{tikz}
\usepackage{bm}
\usepackage{marginnote}
\usepackage{cases}
\usetikzlibrary{shapes,patterns,snakes}
\usetikzlibrary{positioning,fit,calc}

\def\begeq{\begin{equation}}
\def\endeq{\end{equation}}
\def\bc{\begin{center}}
\def\ec{\end{center}}

\title{Why strengthening gap junctions may hinder action potential propagation}

\author{Erin Munro Krull\footnotemark[1] ~\!\! and Christoph B\"{o}rgers\footnotemark[2]}

\begin{document}
\maketitle

\renewcommand{\thefootnote}{\fnsymbol{footnote}}
\footnotetext[1]{Department of Mathematics and Computer Science, Ripon College, Ripon, WI 54971. email: munrokrulle@ripon.edu}
\footnotetext[2]{Department of Mathematics, Tufts University,  Medford, MA 02155. 
email: cborgers@tufts.edu}
\renewcommand{\thefootnote}{\arabic{footnote}}

\begin{abstract} 
Gap junctions are channels in cell membranes allowing ions to pass directly between cells. They connect 
cells throughout the body, including heart myocytes, neurons,  and astrocytes. Propagation mediated by gap junctions can be {\em passive} or {\em active}. In passive propagation, the membrane potential of one cell influences that of neighboring cells without triggering action potentials (APs).  In active propagation, an AP in one cell triggers APs in neighboring cells; this occurs in cardiac tissue and throughout the nervous system. It is known experimentally that there is an ideal gap junction conductance for AP propagation --- weaker {\em or} stronger conductance can block propagation.  We present a theory explaining this phenomenon by analyzing an idealized model that focuses exclusively on  gap junctional and spike-generating currents. We also find a novel type of behavior that we call {\em semi-active} propagation, in which cells in the network are not excitable at rest, but still propagate action potentials.
\end{abstract}



\pagestyle{myheadings}
\thispagestyle{plain}
\markboth{E.\ Munro Krull and C.\ B\"orgers}{Action potential propagation through gap junctions}

\section{Introduction}

Gap junctions are found throughout the body \cite{Bennett:1997vm,Goodenough:2009hr, Rackauskas:2010ts}. In particular, gap junctions are known to connect astrocytes \cite{Bennett:2003bj}, neurons \cite{Bennett:2004uy, Connors:2004jl, Hormuzdi:2004iw, Sohl:2005br} especially during development \cite{Montoro:2004cs, Sutor:2005fi}, and heart myocytes \cite{Bernstein:2006, Kleber:2021, Roerig:2000tj}. More recently, gap junctions were observed in cancer cells \cite{Bonacquisti:2019, Osswald:2015hn}, and between soma and germline cells \cite{Landschaft:2020}.

Gap junctions allow ions to pass between cells, and so the voltage in one cell directly influences the voltage in the neighboring cell. Suppose two excitable cells are connected, a ``trigger" cell and a downstream neighbor, where the voltage of both cells is the resting potential. If there is a small deflection in the trigger cell, then we will see a proportional deflection in the neighbor. This proportion is often referred to as the coupling coefficient \cite{Welzel:2019}. If the gap junction connection is weak, then an action potential (AP) in the trigger cell may only yield a proportional response in the neighbor. This response is called passive propagation, since the neighbor's firing currents are never fully activated. On the other hand, if the gap junction connection is strong, then an AP in the trigger cell may yield an AP in the neighbor. We call this active propagation, since the neighbor's firing currents amplify the response.

Active propagation of APs across gap junctions is seen in cardiac tissue \cite{Bernstein:2006, Kleber:2021, Roerig:2000tj}, and throughout the nervous system \cite{Chorev:2012, Draguhn:1998fm, Mercer:2006bd, Papasavvas:2020, Sheffield:2011bi, Shimizu:2013, Smedowski:2020, Wang:2010hk}. More recently, gap junctions were found  to mediate propagating calcium waves through smooth muscle tissue \cite{Borysova:2018}, neural progenitor cells \cite{Malmersjo:2014}, and possibly astrocytes \cite{Goldberg:2010, Spray:2019, Verhoog:2020} and retinal cells \cite{Kahne:2019}.

In several studies of the effect of gap junction strength on propagation, researchers noted that APs may propagate more easily across a given network if the gap junction strength is set to a medium amount. 
For example, previous experiments in heart tissue show that APs propagate faster when gap junction conductance isn't too weak or too strong \cite{Dhein:2006ej, Kleber:2021, Morley:2005bg, Rohr:1997vt, Rohr:2004gi}. These results are corroborated in heart tissue models, which show 
that there is an ideal gap junction conductance that maximizes the number of downstream neighbors that APs may propagate across \cite{Hubbard:2010jk, Rudy:1987tk, Shaw:1997ws}. Likewise, several neural models involving gap junctions show that APs can propagate when connection strength is set to a medium amount, and AP propagation is blocked when conductance is too high or too low \cite{Gansert:2007gm, Gonzalez-Ramirez:2019, Nadim:2006eq}. To understand these results, we think about increasing gap junction conductance in a network of cells with possibly many downstream neighbors. If the gap junction conductance between cells is zero, clearly APs cannot propagate through the network. As the gap junction conductance is raised, propagation becomes possible. However, if the conductance is raised too much, current shunted to downstream neighbors may prevent APs from propagating.

We describe a theory for understanding when APs can propagate through a network connected by gap junctions. Our theory uses a simplification of the Fitzhugh-Nagumo equations, but we verify that our conclusions hold qualitatively for more realistic biophysical models as well. We also find a novel behavior that we call semi-active propagation, in which cells in the network are not excitable at rest, but still propagate action potentials.

\section{Propagation into a single cell} \label{sec:single_cell}

\subsection{Model of a single cell}

\begin{figure}[ht!]
\centering
\includegraphics[scale=0.35]{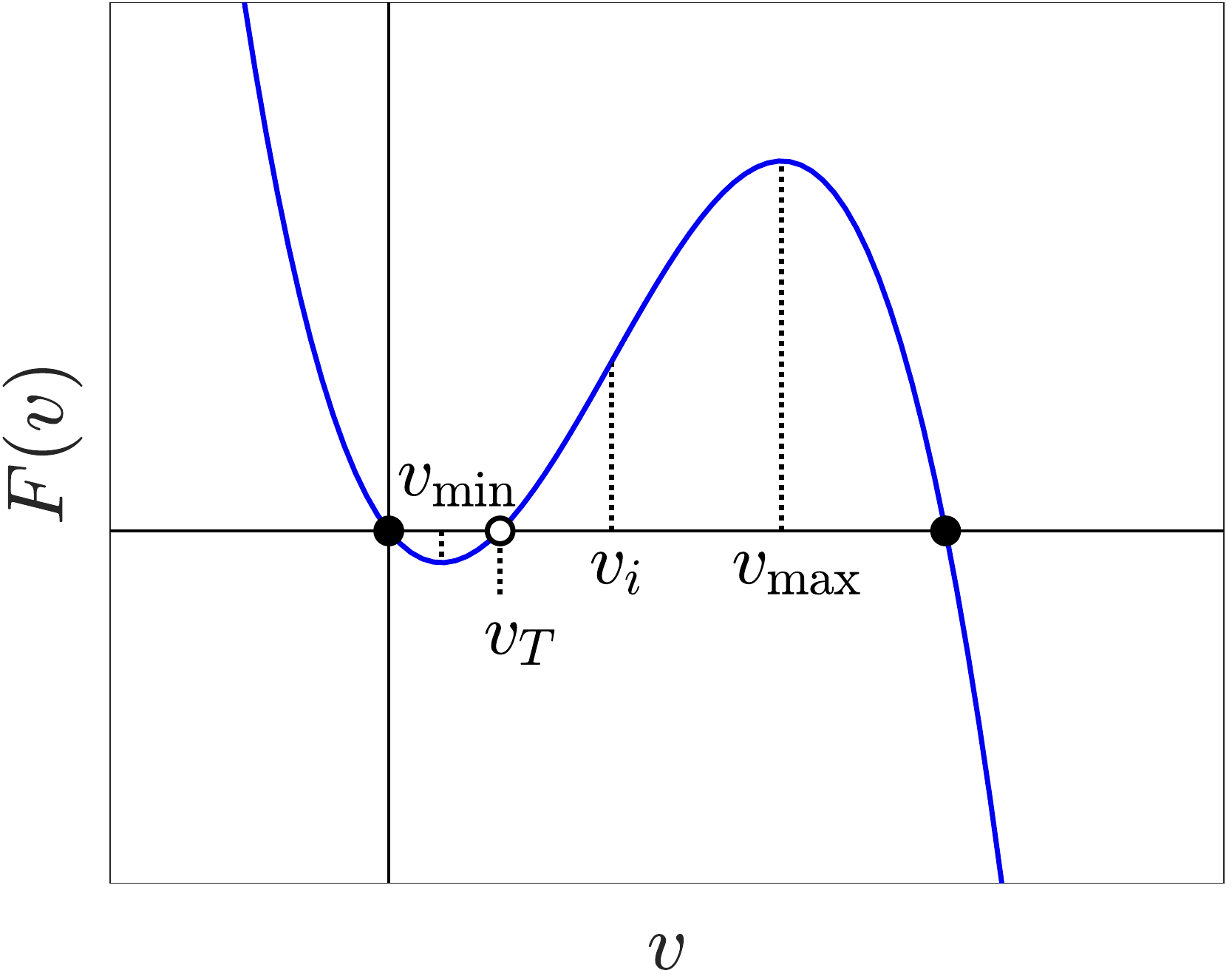}
\caption{The graph of $F$, and some important points.}
\label{fig:VMINVIVMAX}
\end{figure}

We model a single excitable cell using the equation
\begin{equation}
\label{eq:isolated_cell}
\frac{dv}{dt} = v(v-v_T) (1-v) 
\end{equation}
where $v_T$ is a parameter with $0 < v_T < 1$.
We denote the right-hand side of eq.\ (\ref{eq:isolated_cell}) by $F(v)$; see Fig. \ref{fig:VMINVIVMAX}. 
We think of $v$ as a non-dimensionalized membrane potential, and will therefore refer to it as {\em voltage}. We also think of $F(v)$ as  a current (more accurately, current divided by capacitance). Equation (\ref{eq:isolated_cell}) is one half of the FitzHugh-Nagumo model \cite{Ermentrout_Terman_book}. We omit the other half, the slow variable that drives $v$ down when it is high. Our focus, for now, is on the {\em initiation}, not on the {\em termination} of spikes.

Equation (\ref{eq:isolated_cell})  has a stable equilibrium at $v=0$, an unstable one at $v=v_T$, and a stable one at $v=v_F=1$. If $v$ is perturbed slightly from $0$, it returns to $0$, but if it is raised above $v_T$, it rises further to $1$ instead.  In this paper, {\em excitability} always means the co-existence of a stable resting state ($v=0$ in eq.\ (\ref{eq:isolated_cell})) and an unstable {\em threshold} ($v=v_T$ in eq.\ (\ref{eq:isolated_cell})), with the property that perturbation away from the stable resting state and past the threshold results in a large excursion (convergence to the {\em peak voltage} $v=1$ in eq.\ (\ref{eq:isolated_cell})).

The graph of $F$ has  a local minimum $v_{\rm min}$, a local maximum $v_{\rm max}$, and an inflection point at $v_i=(v_T+1)/3$. We assume $v_T < v_i$, which is equivalent to 
$$
v_T < \frac{1}{2}.
$$
For neurons and myocytes, this is a natural assumption: The firing threshold is closer to the resting voltage ($v=0$) than to the peak voltage ($v=1$). So 
$$
v_{\rm min} < v_T < v_i < v_{\rm max} 
$$
as illustrated in Fig.\ \ref{fig:VMINVIVMAX}.

\subsection{Gap junctional connections with neighbors held at fixed voltage}

\begin{figure}[ht!]
\centering
\includegraphics[scale=1.0]{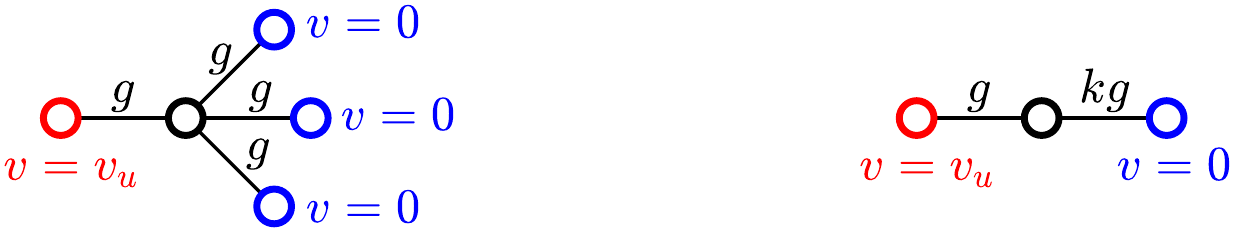}
\caption{Left: Cell with upstream and downstream neighbors at fixed voltage. Right: Downstream neighbors collapsed into 
a single cell.}
\label{fig:NETWORK}
\end{figure}

Imagine that the cell in the previous section has $k+1$ gap junctionally connected neighbors, of which one, referred to as {\em upstream}, has voltage $v_u$, and $k$ others, referred to as {\em downstream}, are at rest with $v=0$.  The upstream and downstream voltages are assumed fixed for now. Only the central cell has dynamics. We assume that all gap junctional conductances have the same strength $g$. See Fig.\ \ref{fig:NETWORK}, left panel, for an illustration with $k=3$.

\begin{figure}[ht!]
\centering
\includegraphics[scale=0.35]{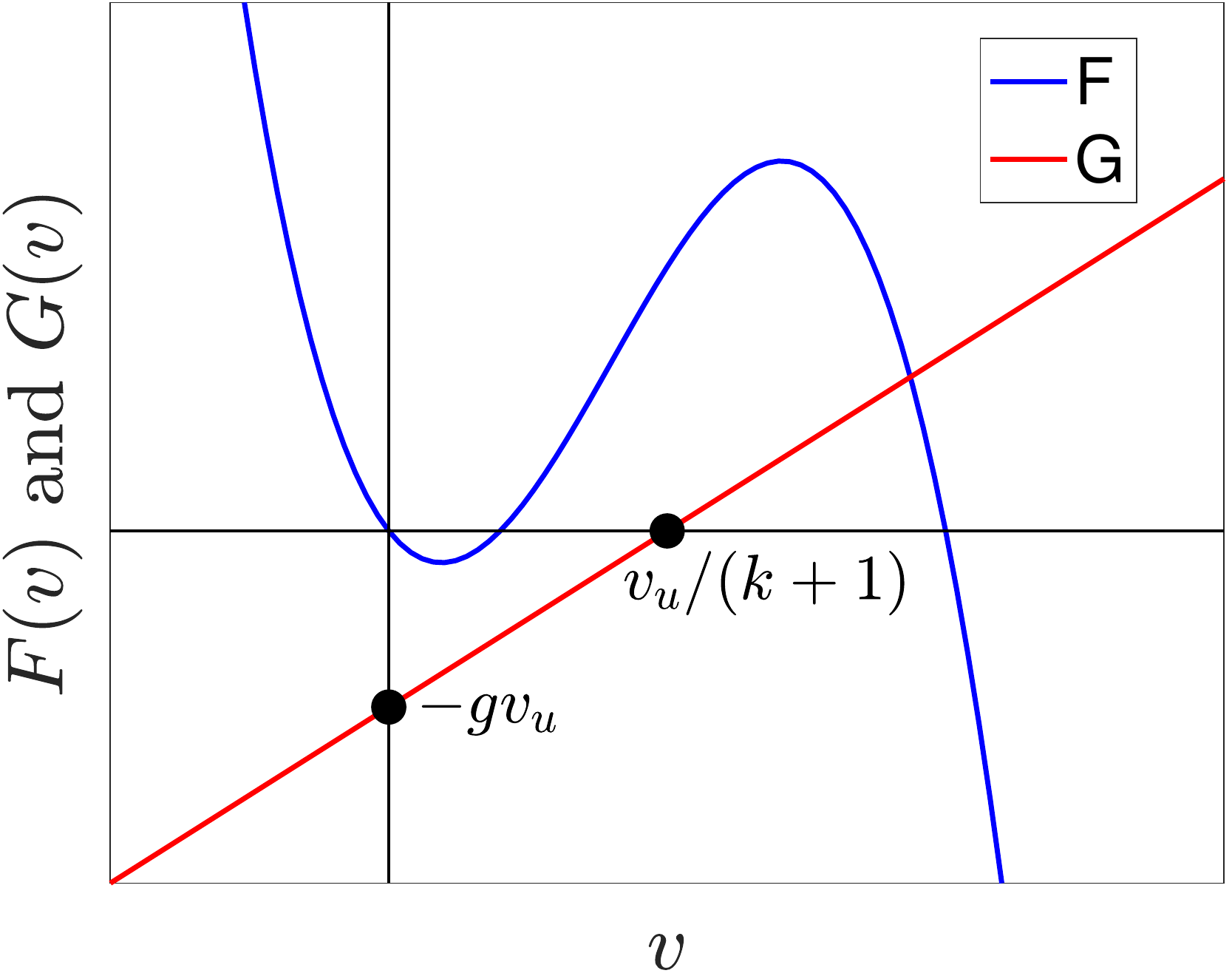}
\caption{$F$  and $G$.}
\label{fig:GRAPH_OF_G}
\end{figure}

Our model now becomes
\begin{equation}
\label{eq:with_gap_junctions}
\frac{dv}{dt} = v(v-v_T) (1-v) + g (v_u-v)  -gkv.
\end{equation} 
We assume here, as is customary, that gap junctions are governed by Ohm's law --- current is proportional 
to voltage difference. The term $g(v_u-v)$ models the effect of the upstream cell, and $-g k v$ models the effect  of the downstream cells. With the notation
\begin{equation}
\label{defG}
G(v) =  g (k+1) v - g v_u, 
\end{equation}
 eq.\ (\ref{eq:with_gap_junctions}) becomes
$$
\frac{dv}{dt} = F(v) - G(v).
$$
Figure \ref{fig:GRAPH_OF_G} shows $F$ and $G$ in one figure. The graph of $G$ intersects the horizontal axis at $(v_u/(k+1),0)$, and the vertical axis at $(0,-g v_u)$.  

In general, there might be more upstream cells than just one. There might also be more or fewer downstream cells, coupled to the central cell with varying conductances. What matters  is that $g$ is the total gap junctional conductance with which the central cell is coupled to upstream  cells, and $k g$ is the total gap junctional conductance with which the central cell is coupled to downstream  cells. Therefore we won't assume $k$ to be an integer any more,  simply thinking of it as the ratio of total downstream conductance over total upstream conductance, and collapse our  network picture as in the right panel of Fig.\ \ref{fig:NETWORK}.

\subsection{All neighbors at rest} 
\label{sec:neighbors_at_rest}

\begin{figure}[ht!]
\centering
\includegraphics[scale=0.84]{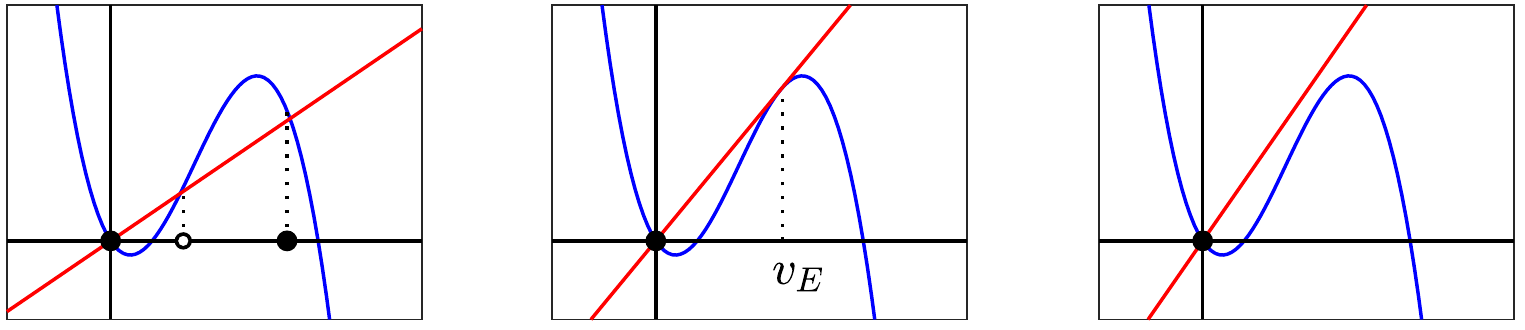}
\caption{$F$  and $G$ for $v_u=0$, different values of $g(k+1)$.}
\label{fig:ZERO_VU_2}
\end{figure}

We first think about the case when the upstream neighbors of the central cell, just like the downstream neighbors, are at rest: $v_u=0$. We will discuss the equilibria of the central cell in that case. The formula for $G(v)$ is now 
$$
G(v) = g (k+1) v.
$$
Fixed points of eq.\ (\ref{eq:with_gap_junctions}) are solutions of 
\begin{equation}
\label{eq:fps}
F(v) = G(v).
\end{equation}
If $g(k+1)$ is not too large, then there are three equilibria, depicted in the left panel of Fig.\ \ref{fig:ZERO_VU_2}. As $g(k+1)$ increases, the graph of $G$ (shown in red) becomes steeper. This raises the threshold and lowers the peak voltage, but doesn't affect the resting voltage. Eventually the threshold and peak voltage collide at a voltage which we denote by $v_E$; see middle panel of Fig.\ \ref{fig:ZERO_VU_2}. For larger $g(k+1)$, there is only one equilibrium at $v=0$. 

We say that the cell is {\em excitable} if there are three equilibria --- an unstable threshold surrounded by two stable equilibria. So gap junctional connections with neighbors held at rest, if they are strong or numerous enough, can make excitability disappear altogether.

\subsection{Firing triggered by raising the upstream voltage}
\label{sec:triggering}

\begin{figure}[ht!]
\centering
\includegraphics[scale=0.75]{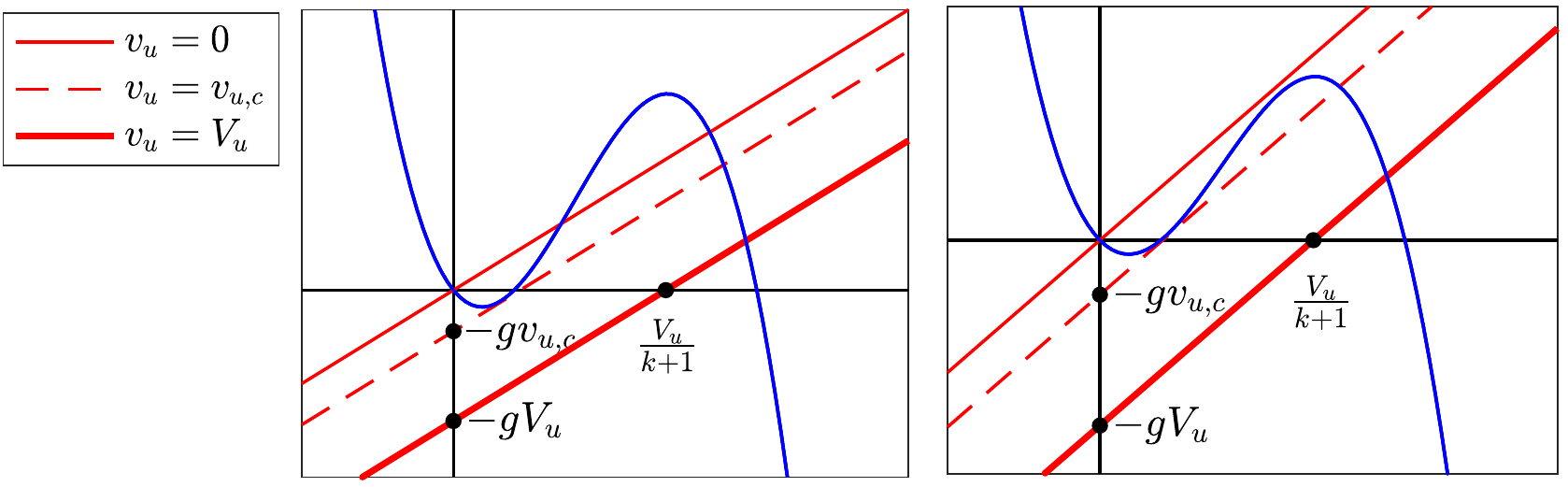}
\caption{Raising $v_u$ shifts the graph of $G$.}
\label{fig:THREE_LINES}
\end{figure}

As $v_u$ is raised from $0$ to some maximum value $V_u>0$, the graph of $G$ moves downward. 
Depending on the values of $V_u$, $k$, and $g$, the threshold and resting voltages may eventually collide and annihilate each other, leaving the peak voltage as the only equilibrium. We denote by $v_{u,c}$ the value of $v_u$ at which the collision of threshold and resting voltage occurs, $0 < v_{u,c} < V_u$; see Fig.\ \ref{fig:THREE_LINES}.

As the right panel of Fig.\ \ref{fig:THREE_LINES} shows, even when the central cell  with all neighbors at rest is not excitable because it is too strongly connected to resting neighbors, raising $v_u$ may restore excitability, by making the threshold and peak voltage equilibria re-appear via a blue sky bifurcation. Then, at $v_u=v_{u,c}$, the resting voltage and threshold collide and disappear.

\vskip 5pt
\begin{definition}
We say that the central cell {\em fires}, or has a {\em firing response}, when $v_u$ rises from $0$ to $V_u$, if there is a threshold-rest collision at some critical value $v_u = v_{u,c}$ with $0 < v_{u,c} < V_u$.
\end{definition}

\vskip 5pt
\begin{definition}
We refer to the segment between $v=v_{\rm min}$ and $v=v_i$ on the graph of $F$ as the {\em critical segment}.
\end{definition}
\vskip 5pt

The critical segment is indicated in the left panel of Fig.\ \ref{fig:GMIN} as a bold blue curve. Any threshold-rest collision must occur on the critical segment. The following proposition follows from Fig.\ \ref{fig:THREE_LINES}. 

\vskip 5pt
\begin{proposition}
\label{prop:collision_conditions}
Let $g>0$ and $k \geq 0$. The central cell fires when $v_u$ is raised from $0$ to $V_u$ if and only if the line $L_{g,k}$ through ($V_u/(k+1),0)$ and $(0, - g V_u)$ satisfies the following two conditions.
\begin{itemize}
\item[(1)] $L_{g,k}$ lies
strictly below the critical segment, and 
\item[(2)] $L_{g,k}$ has slope $<F'(v_i)$.
\end{itemize}
\end{proposition}

\vskip 5pt
Note that $L_{g,k}$ is the graph of $G$ when $v_u=V_u$. We will call $L_{g,k}$ {\em admissible} if it satisfies the conditions in Proposition \ref{prop:collision_conditions}. The left panel of Fig.\ \ref{fig:GMIN} shows an example.

\subsection{Firing depends non-monotonically on $g$}\label{sec:non-mon}

We use  Proposition \ref{prop:collision_conditions} to derive a detailed description of  the set of the parameter pairs $(g,k)$ for which there is a firing response of the central cell. First, we note that $L_{g,k}$ cannot be admissible if $V_u \leq v_T$. We will therefore assume  
$$
V_u > v_T.
$$

\begin{figure}[ht!]
\centering
\includegraphics[scale=0.8]{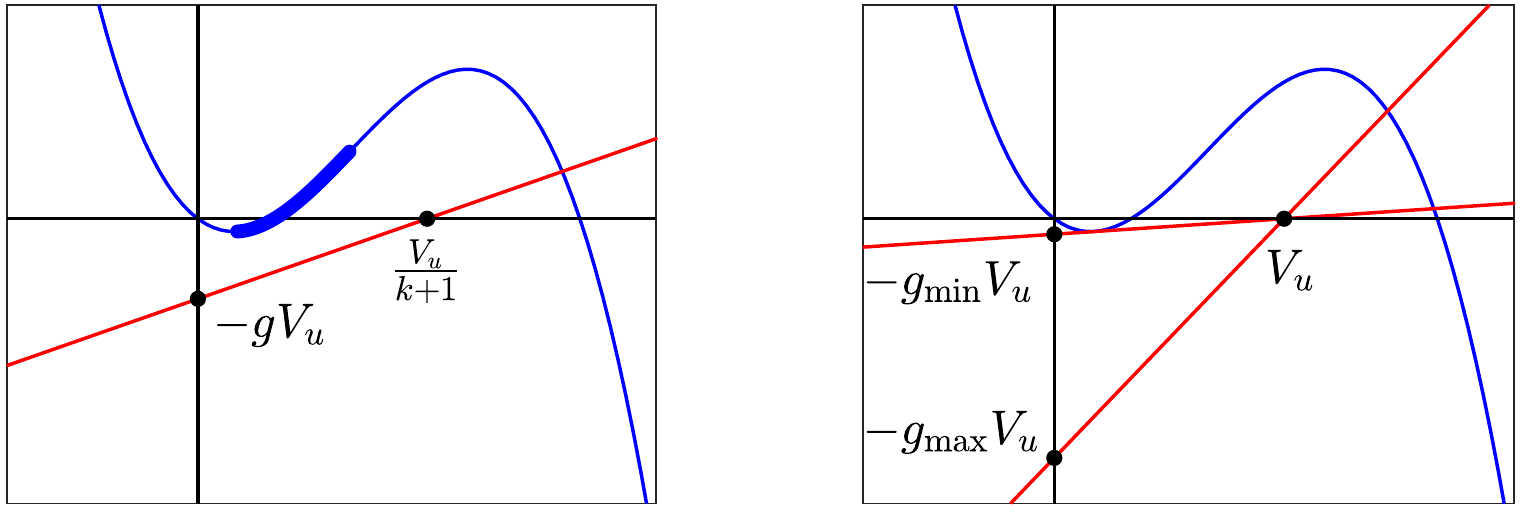}
\caption{An admissible line $L_{g,k}$ (left), and the minimal and maximal values of $g$ allowing signal propagation (right). The bold blue curve segment in the left panel is the {\em critical segment}. }
\label{fig:GMIN}
\end{figure}

Furthermore, since lowering $k$ rotates $L_{g,k}$ clockwise around the point $(0,-gV_u)$, lowering $k$ makes it easier to fire: For $L_{g,k}$ to be admissible for a given $g$, $L_{g,0}$ has to be admissible. (This makes sense, since lowering $k$ lowers the current shunted to the downstream neighbors.) Therefore we first analyze the case $k=0$. Note that $L_{g,0}$ passes through $(V_u,0)$. For $L_{g,0}$ to be admissible, $g$ must satisfy 
$$
g_{\rm min} < g < g_{\rm max},
$$
where the definitions of $g_{\rm min}$ and $g_{\rm max}$ are illustrated by the right panel of Fig.\ \ref{fig:GMIN}.  The line $L_{g_{\rm max},0}$ (the steeper of the two red lines in the figure) has slope $F'(v_i)$, so $g_{\rm max} = F'(v_i)$. The line $L_{g_{\rm min},0}$ is tangent to the graph of $F$.

Now we consider $k>0$. Suppose $g$ lies between $g_{\rm min}$ and $g_{\rm max}$. Then $L_{g,0}$ is admissible. Other admissible lines, for the same value of $g$, are obtained by rotating this line around the point $(0, -g V_u)$ in the counter-clockwise direction until it touches the critical segment or reaches slope $F'(v_i)$ --- whichever comes first. 

\begin{figure}[ht!]
\centering
\includegraphics[scale=0.8]{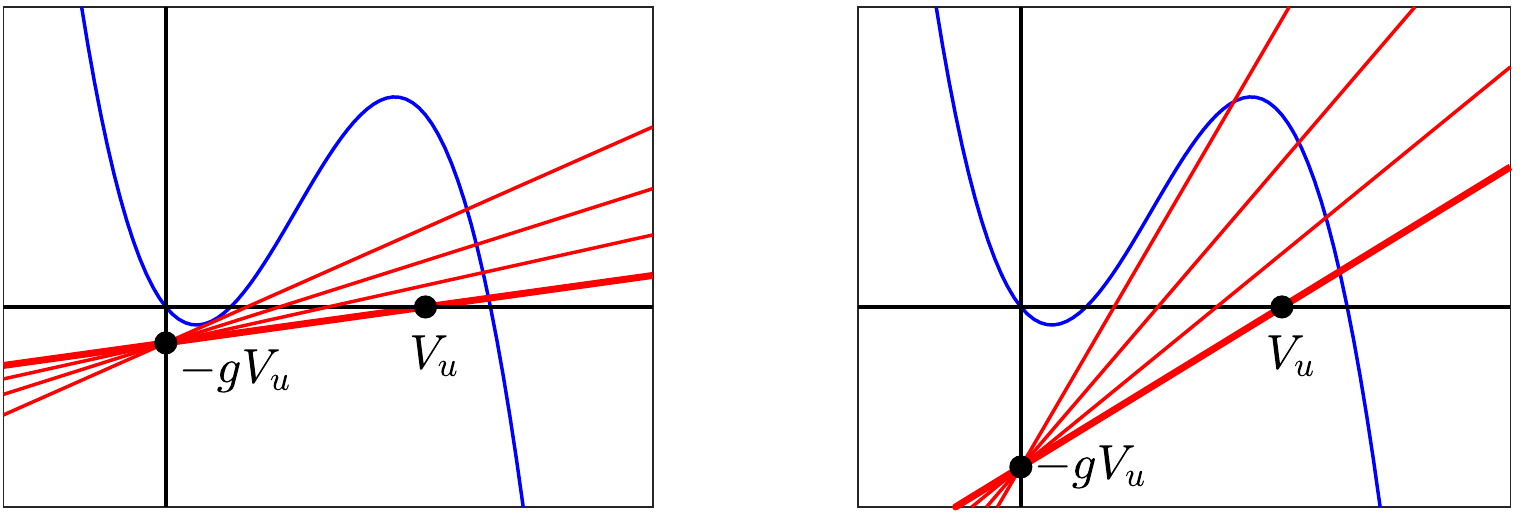}
\caption{Left: An example where the maximum $k$ occurs when $L_{g,k}$ is tangent to the graph of $F$. Right: An example where the maximum $k$ occurs when $L_{g,k}$ reaches $F'(v_i)$.}
\label{fig:WHICHEVER_COMES_FIRST}
\end{figure}

Which of these two does come first depends on $g$. If $g$ is smaller than some threshold value $g_\ast$ (discussed further below), then as the line rotates counter-clockwise around $(0,-g V_u)$, it  touches the critical segment before it reaches slope $F'(v_i)$; see the left panel of Fig.\ \ref{fig:WHICHEVER_COMES_FIRST}.  If $g$ is larger than $g_\ast$, the line reaches slope $F'(v_i)$ before it touches the critical segment; see the right panel of the figure.

From Fig.\ \ref{fig:WHICHEVER_COMES_FIRST}, we see that $g_\ast$ is determined by the tangent to the cubic at $(v_i,F(v_i))$, which intersects the vertical axis at  $-g_\ast V_u$. That is, 
\begin{equation}
\label{def_g_ast}
g_\ast =  \frac{F'(v_i) v_i - F(v_i)}{V_u}.
\end{equation}

Is it  possible that $g_\ast \geq g_{\rm max} = F'(v_i)$? According to (\ref{def_g_ast}), this means
$$
\frac{F'(v_i) v_i - F(v_i)}{V_u} \geq F'(v_i)
$$
or 
\begin{equation}
\label{eq:V_u_klitzeklein}
V_u \leq  v_i - \frac{F(v_i)}{F'(v_i)}.
\end{equation}
The right-hand side of this inequality is the $v$-intercept of the tangent to the cubic at $(v_i,F(v_i)))$. So $g_\ast \geq g_{\rm max}$ only if $V_u$ happens to be very small, just above $v_T$.

For every $g \in (g_{\rm min},g_{\rm max})$, the slope 
$$m=g(k+1)$$ 
of the admissible line $L_{g,k}$ passing through $(0,-g V_u)$ and $(V_u/(k+1),0)$  is bounded above, and the upper bound on $m$ translates into  an upper bound on $k$.  So the parameter regime in which there is a firing response is of the form 
 $$
g_{\rm min}<g<g_{\rm max}, ~~~~ 0 \leq k < k_{\max}(g),
 $$
 and our goal will be to understand what the graph of the function $k_{\rm max}$ looks like. 
 
When $g_\ast \leq g < g_{\rm max}$, the constraint on $m$ is 
$m < F'(v_i) = g_{\rm max}$, and therefore 
\begin{equation}
\label{eq:kmaxformula}
k_{\rm max} (g) =  \frac{F'(v_i)}{g} - 1.
\end{equation}
For $g_{\rm min} < g < g_\ast$, $k_{\rm max}(g)$ is obtained by computing the tangent to the critical segment passing through $(0,-g V_u)$. The two pieces match up continuously at $g_\ast$, and we set 
 $$k_\ast = k_{\rm max}(g_\ast).$$

\begin{figure}[ht!]
\centering
\includegraphics[scale=.7]{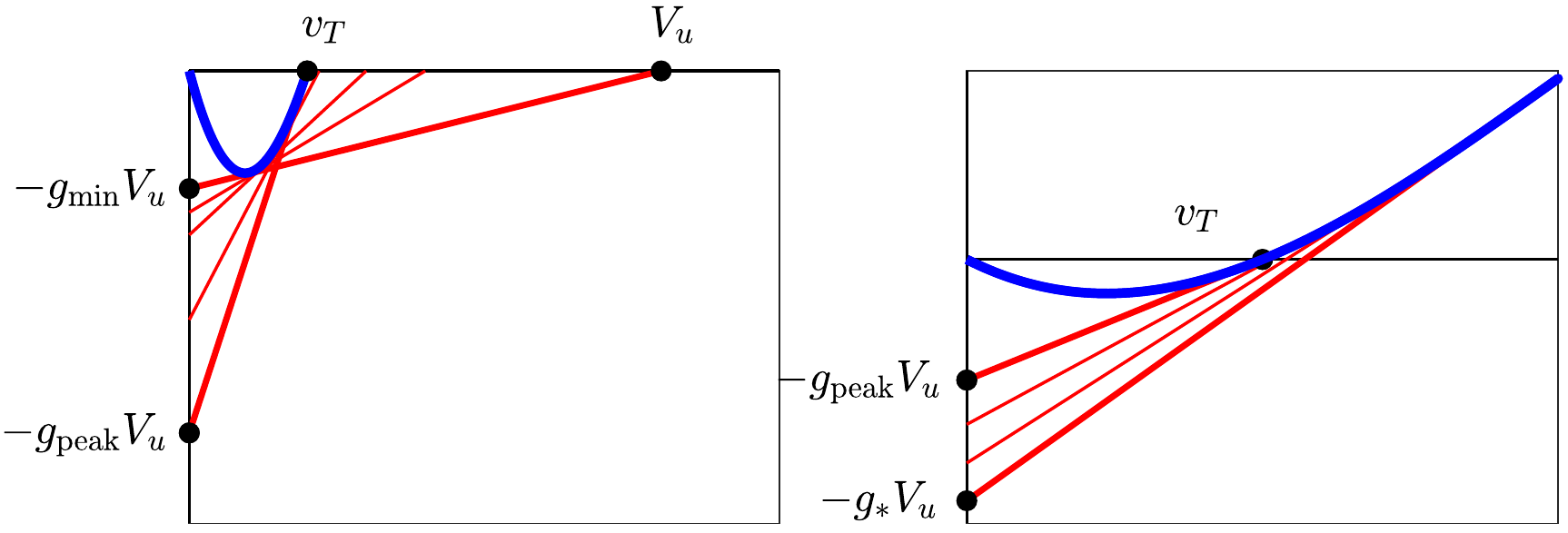}
\caption{As $g$ increases from $g_{\rm min}$ to $g_\ast$, the $v$-intercept of the tangent through $(0,-g V_u)$ first moves left (left panel), then right (right panel). Therefore $k_{\rm max}(g)$ first increases, then decreases.} 
\label{fig:RISES_FALLS}
\end{figure} 

\begin{figure}[ht!]
\centering
\includegraphics[width=5in]{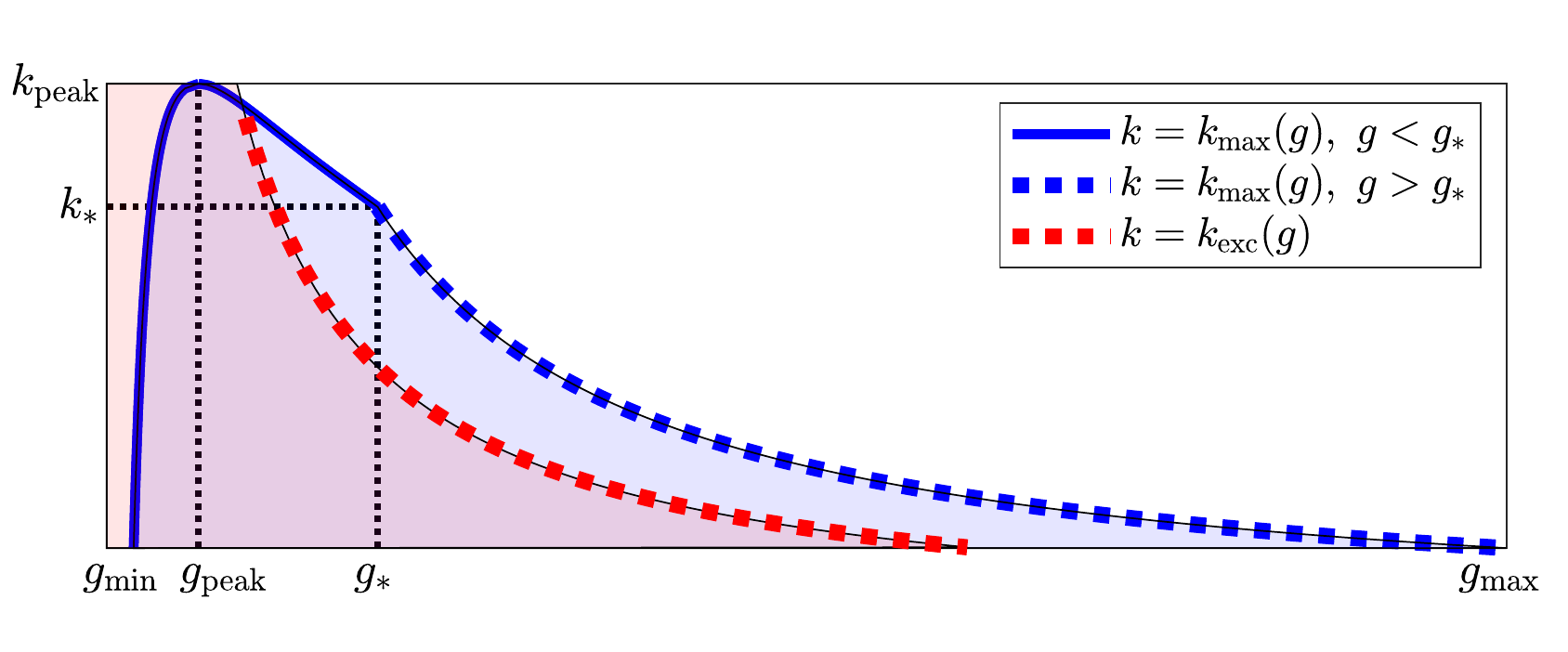}
\caption{The region in the $(g,k)$-plane in which firing occurs (blue), and the region in which the cell is excitable with all neighbors at rest (red), for $v_T=0.15$, $V_u=1$. 
}
\label{fig:G_K}
\end{figure}

We will now follow what happens to the $v$-intercept, $V_u/(k_{\rm max}(g) +1)$, as $g$ increases. As shown in Fig.\ \ref{fig:RISES_FALLS}, initially the $v$-intercept falls as $g$ increases, until it reaches $v_T$. As $g$ continues to increase, the $v$-intercept rises again. Therefore $k_{\rm max}$ rises first, then falls. The leftmost $v$-intercept, $v_T$, occurs for a value of $g$ that we call $g_{\rm peak}$, and $k_{\rm max}$ reaches its maximum at $k_{\rm peak} = k_{\rm max}(g_{\rm peak})$. We have $V_u/(k_{\rm peak}+1) = v_T$ and therefore
\begin{equation}
\label{eq:k_peak}
k_{\rm peak} = \frac{V_u}{v_T} - 1.
\end{equation}
and since $g_{\rm peak} (k_{\rm peak}+1) = F'(v_T)$, 
\begin{equation}
\label{peak_equation}
g_{\rm peak} =  \frac{F'(v_T)}{k_{\rm peak}+1} = \frac{F'(v_T) v_T}{V_u}.
\end{equation}
Our conclusions are summarized in Fig.\ \ref{fig:G_K}. When $V_u$ is so small that $g_\ast \geq g_{\rm max}$, the dashed part of the blue boundary in Fig.\ \ref{fig:G_K} is simply absent.

We also indicate in Fig.\ \ref{fig:G_K} the parameter pairs $(g,k)$ for which the central cell is excitable
when the upstream cell is at rest, $v_u=0$. The discussion in Section \ref{sec:neighbors_at_rest} shows that this region is given by $g(k+1)<F'(v_E)$, that is, 
\begin{equation}
\label{excitability_boundary}
k < \frac{F'(v_E)}{g} - 1 . 
\end{equation}
We denote the right-hand side of (\ref{excitability_boundary}) by $k_{\rm exc}(g)$. 

\subsection{Active, semi-active, and passive propagation}

We call the region of parameter pairs $(g,k)$ with 
$$
g_{\rm min} < g < g_{\rm max}, ~~~ 0 \leq k < k_{\rm max}(g), ~~~ k < k_{\rm exc}(g)
$$
the region of {\em active propagation}. It appears in purple in Fig.\ \ref{fig:G_K}. When $(g,k)$ lies in this
region, the central cell is excitable for $v_u=0$, and fires when $v_u$ is set to $V_u$. Further, we call the region of pairs $(g,k)$ 
$$
g_{\rm min} < g < g_{\rm max}, ~~~ 0 \leq k < k_{\rm max}(g), ~~~ k  \geq k_{\rm exc}(g)
$$
the region of {\em semi-active propagation}.  When $(g,k)$ lies in this region, the central cell is not 
excitable for $v_u=0$, but as $v_u$ rises from $0$ to $V_u$, first excitability is restored, and then there is a firing response. Pairs $(g,k)$ that don't trigger firing of the central cell when $v_u$ is set to $V_u$ are said to lie in the region of {\em passive propagation}. Raising $v_u$ from $0$ to $V_u$ still causes the voltage of the central cell to rise, but does not trigger firing.

\begin{figure}[ht!]
\centering
\includegraphics[scale=0.6]{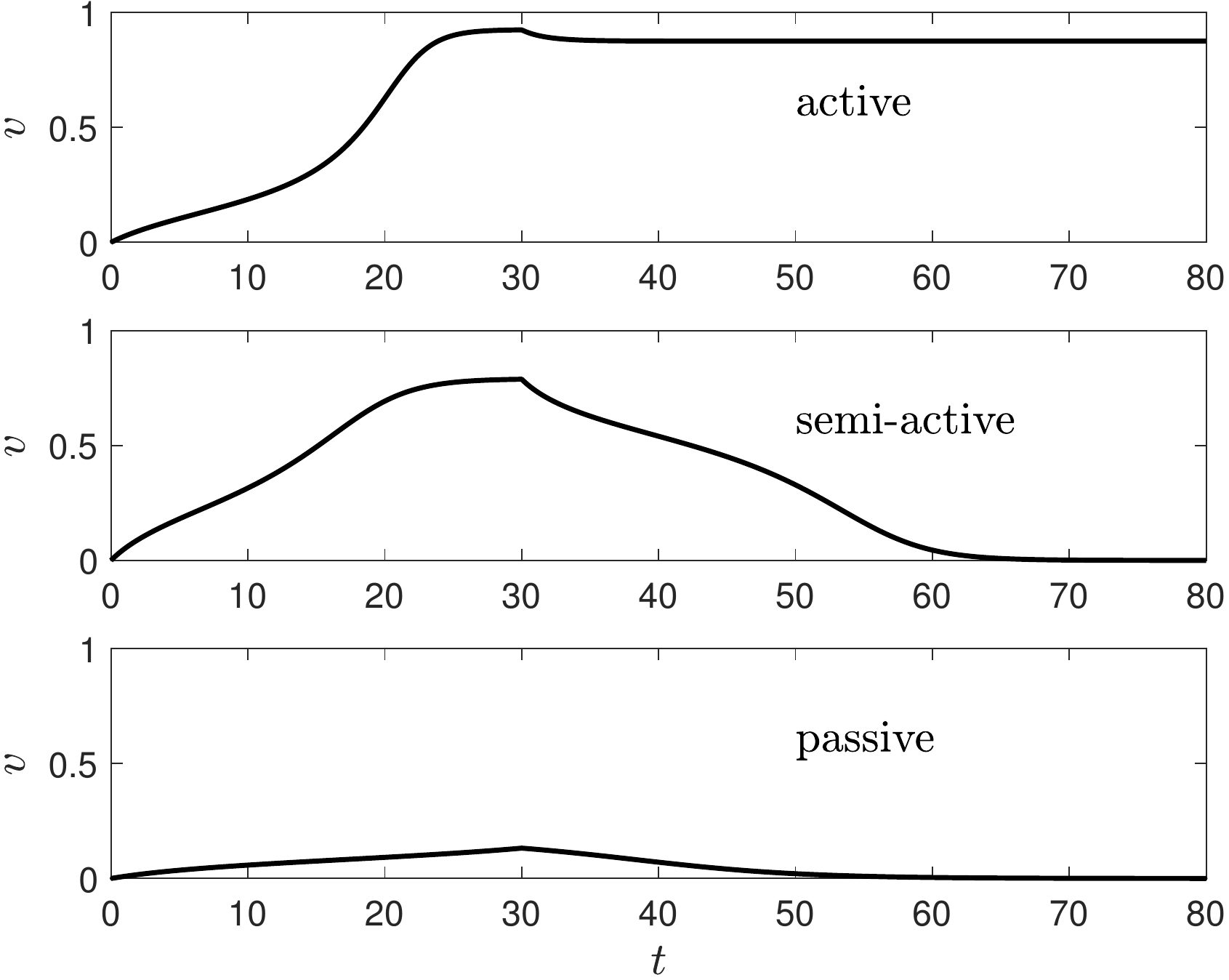}
\caption{Responses to setting $v_u=1$ for $0 \leq t \leq 30$,  then setting it back to $0$, for $k=2$ and $g=0.03$ (top), $g=0.07$ (middle), and $g=0.01$ (bottom).}
\label{fig:THREE_REGIMES}
\end{figure}

To illustrate the difference between the three parameter regimes, we show in Fig.\ \ref{fig:THREE_REGIMES} solutions of (\ref{eq:with_gap_junctions}) with $v(0)=0$ and 
$$
v_u = \left\{ \begin{array}{cl} 1 & \mbox{for $0 \leq t \leq 30$} \\
0 & \mbox{for $t > 30$}.
\end{array}
\right.
$$
In the regime of active propagation, the voltage of the central cell rises, but does not return to $0$ when $v_u$ is set back to zero. (Notice that our model includes no explicit mechanism for bringing the voltage back down after firing.) In the regime of semi-active propagation, the voltage of the central cell rises nearly as much, but eventually returns to zero when $v_u$ is set back to zero, since for $v_u=0$, there is no non-zero equilbirium --- the central cell is not excitable. In the regime of passive propagation, the voltage rises much less, and returns to zero soon after $v_u$ is  switched  back to zero. 

\subsection{Spike height can depend discontinuously on g and k}
\label{subsec:disc}

For $g>0$, and $k \geq 0$, let $v=v(t)$ be the solution of (\ref{eq:with_gap_junctions}) with $v_u=V_u$ and 
$v(0)=0$. Let
\begin{equation}
\label{eq:def_v_infty}
v_\infty = \lim_{t \rightarrow \infty} v(t).
\end{equation}
In other words, $v_\infty$ is the voltage obtained by setting the voltage in the upstream cell to $V_u$, then letting the voltage in the central cell equilibrate, starting at $v=0$; so $v_{\infty}$ is  the smallest equilibrium. We view $v_\infty$ as a function of $g$ and $k$.

\vskip 5pt
\begin{theorem} $v_\infty(g,k)$ has jump discontinuities at parameter pairs $(g,k)$ with $g_{\rm min} \leq g < g_\ast$ and $k=k_{\rm max}(g)$, and is continuous everywhere else. That is, $v_\infty$ is discontinuous along the solid blue curve in Fig.\ \ref{fig:G_K}, but nowhere else.
\end{theorem}

\vskip 5pt
\begin{proof} Let $g>0$ and $k \geq 0$. As before, denote by $L_{g,k}$ the straight line through $(0,-g V_u)$ and $(V_u/(k+1),0)$. Then $v_{\infty}(g,k)$ is the $v$-coordinate of the left-most intersection point of $L_{g,k}$ and the graph of $F$. From Fig.\ \ref{fig:THREE_LINES} we see that $v_\infty$ depends discontinuously on $(g,k)$ if and only if $L_{g,k}$ is tangent to the graph of $F$ along the critical segment. The $(g,k)$ for which $L_{g,k}$ is tangent to the graph of $F$ at the critical segment are precisely the ones that lie on the solid part of the blue curve in Fig.\ \ref{fig:G_K}.
\end{proof}

\begin{figure}[ht!]
\centering
\includegraphics[scale=0.75]{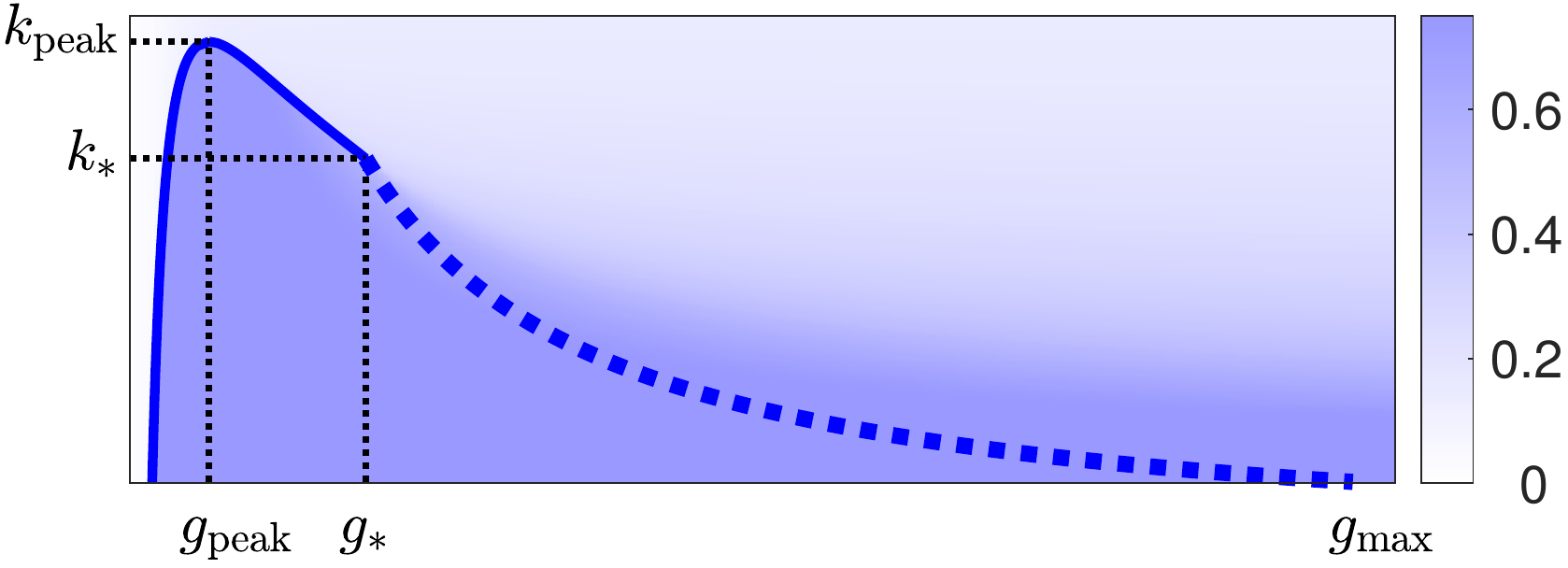}
\caption{Heat plot of $v_\infty$ as a function of $g$ and $k$ for $v_T=0.15$, $V_u=1$. The blue curves are from Fig.\ \ref{fig:G_K}. We see that $v_\infty$ is discontinuous along the solid blue curve, but not along the dashed blue curve.}
\label{fig:HEAT_PLOT}
\end{figure}

\vskip 5pt
Figure \ref{fig:HEAT_PLOT} illustrates this result. 

\begin{figure}[ht!]
\centering
\hskip -30pt
\includegraphics[scale=.55]{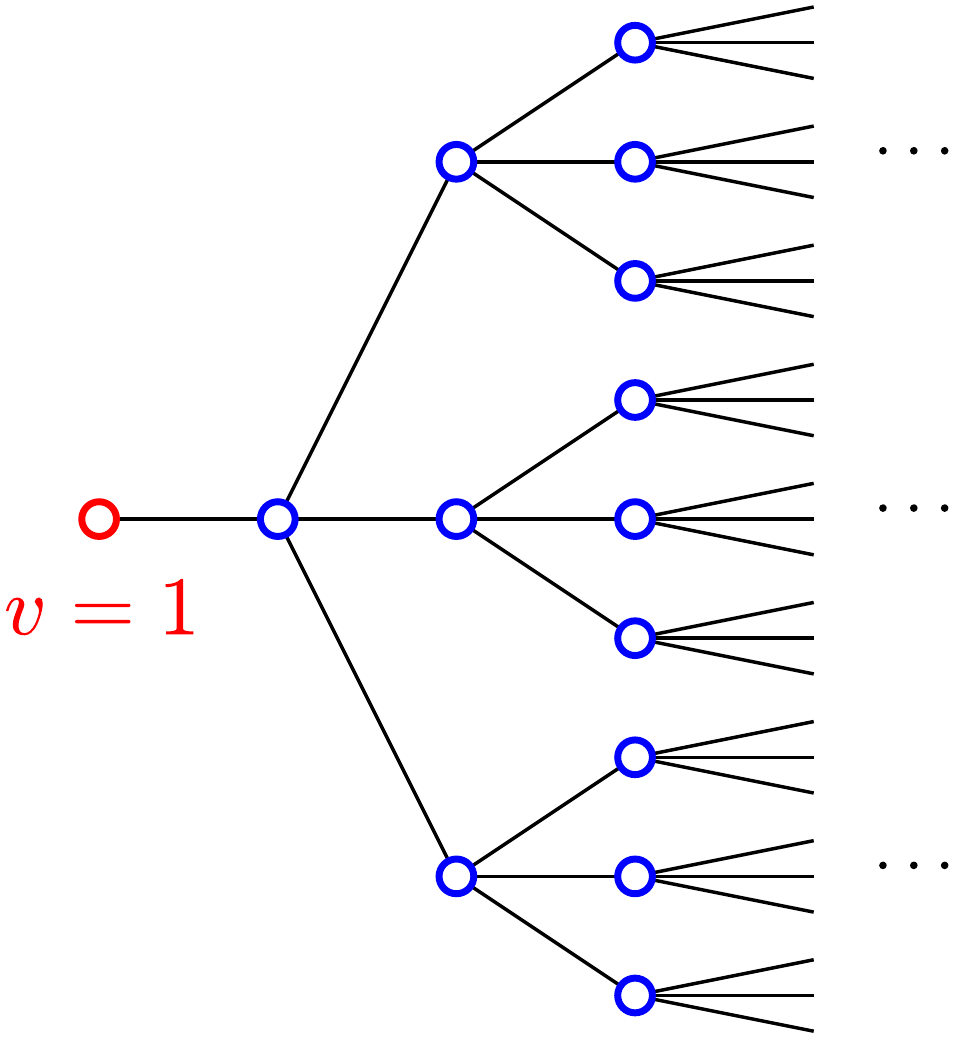}  
\hskip 130pt ~
\vskip -100pt
\hskip 220pt 
\includegraphics[scale=.6]{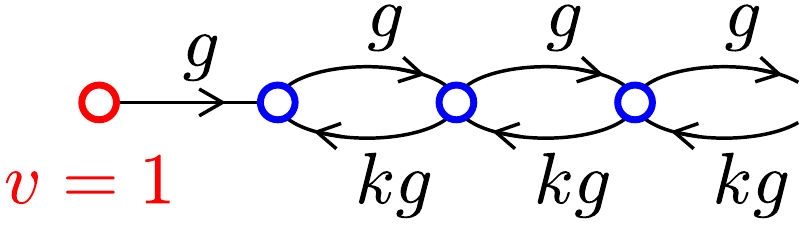}
\vskip 80pt
\caption{Left: tree network. Each cell has one upstream neighbor and $k$ downstream neighbors. Here $k=3$. Propagation is triggered by setting the voltage in the left-most cell, highlighted in red, to $v=1$. Right: equivalent collapsed network.}
\label{fig:TREE}
\end{figure}

\pagebreak
\section{Propagation through a tree} \label{sec:tree}

\subsection{Model}

Now we think about propagation through a tree of  cells connected by gap junctions, as shown in Fig.\ \ref{fig:TREE}. Assume that the cells are identical, and that each cell has exactly one upstream and $k$ downstream neighbors. Propagation is triggered by raising the voltage of the leftmost cell to $v=1$. 

By symmetry, the voltages in cells that are vertically aligned will  be identical. We therefore collapse vertically aligned cells into a single cell, and obtain the equivalent chain shown in the right panel of Fig.\ \ref{fig:TREE}. In this chain, gap junction conductances are not symmetric: $g$ for input from upstream, and $kg$ for input from downstream. As before, we now drop the assumption that $k$ is an integer, viewing it instead as a measure of how densely connected the network is. 

We make the additional simplification that each cell responds only to the voltage in its upstream neighbor and behaves as though its downstream neighbors are held at zero voltage; its voltage therefore rises to $v_\infty$, the smallest equilibrium point of eq.\ \ref{eq:with_gap_junctions}. Propagation through the chain amounts to iterating the map 
$$
\varphi: ~~ v_u \mapsto v_\infty,
$$
which maps $[0,1]$ into $[0,1]$. 

\subsection{Iterating the map $\mathbf{\varphi}$} 

\begin{figure}[ht!]
\centering
\includegraphics[scale=0.6]{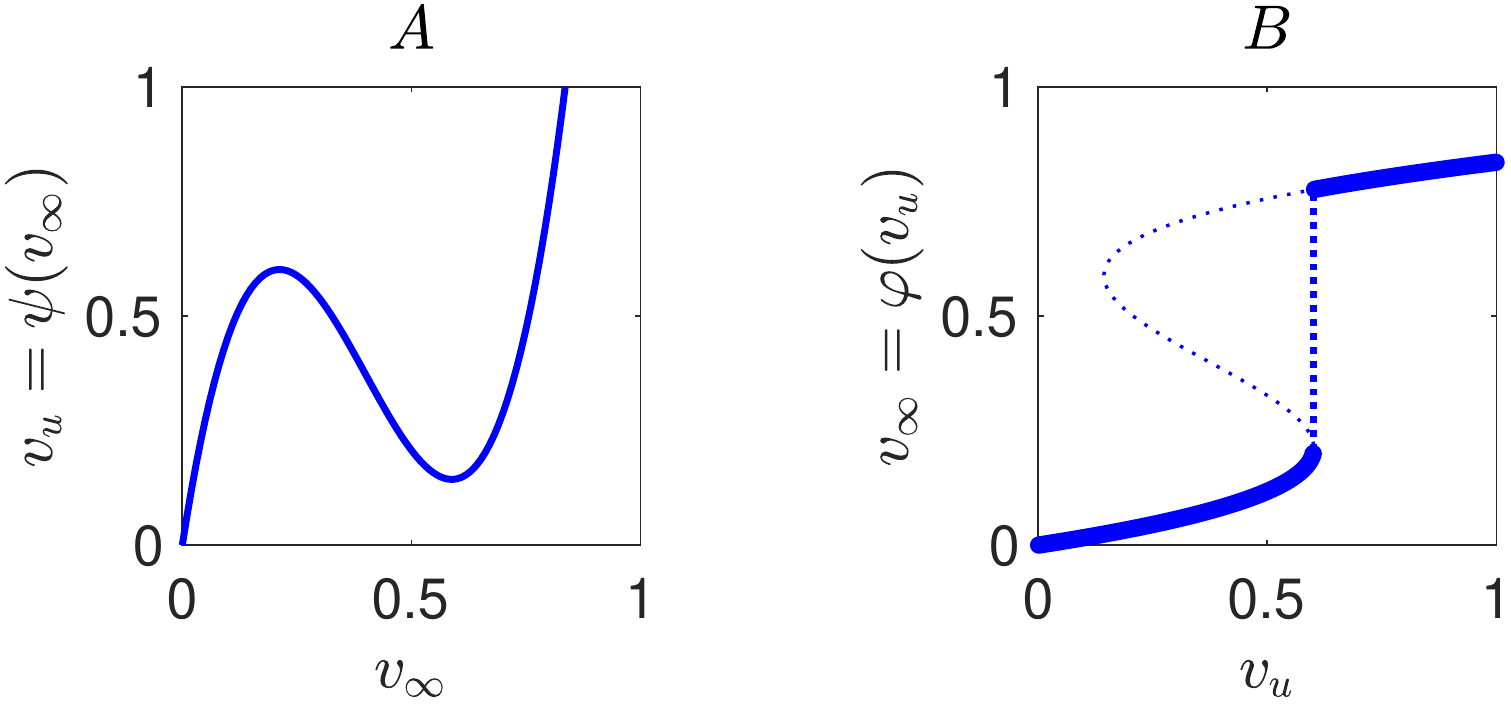}
\caption{A: The function $\psi$ mapping equilibria of eq.\ (\ref{eq:with_gap_junctions}) onto $v_u$. B (bold): The function $\varphi$ mapping $v_u$ onto the smallest equilibrium of eq.\ (\ref{eq:with_gap_junctions}).}
\label{fig:PLOT_PHI_PSI}
\end{figure}

A slight complication in our analysis is that the map $\varphi$ need not be continuous. To see this, first note that given any equilibrium $v_\infty$ of eq.\ (\ref{eq:with_gap_junctions}), we can calculate $v_u$ easily: 
$$
v_u =  - \frac{F(v_\infty)}{g} + (k+1) v_\infty.
$$
We denote the right-hand side of this equation by $\psi(v_\infty)$; Fig.\ \ref{fig:PLOT_PHI_PSI}A shows an example. The graph of $\varphi$ is obtained by swapping the $v_\infty$- and $v_u$-axes, and using the smaller of the two possible values of $v_\infty$ for a given $v_u$ when there is ambiguity; see Fig.\ \ref{fig:PLOT_PHI_PSI}B.

\begin{figure}[ht!]
\centering
\includegraphics[scale=0.75]{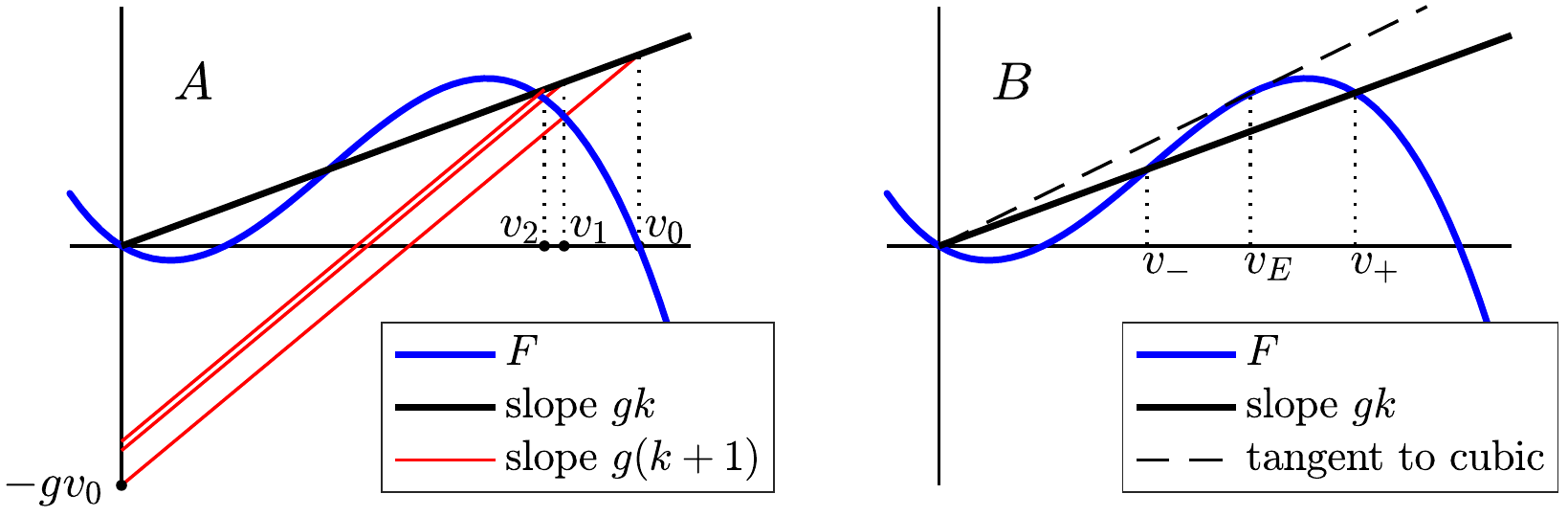}
\caption{A: Construction of $v_1 = \varphi(v_0)$ and $v_2 = \varphi(v_1)$, starting with $v_0=1$. \\ B: Definitions of $v_E$, $v_+$, $v_-$.}
\label{fig:CONSTRUCTING_PHI}
\end{figure}

\begin{figure}[ht!]
\centering
\includegraphics[scale=0.78]{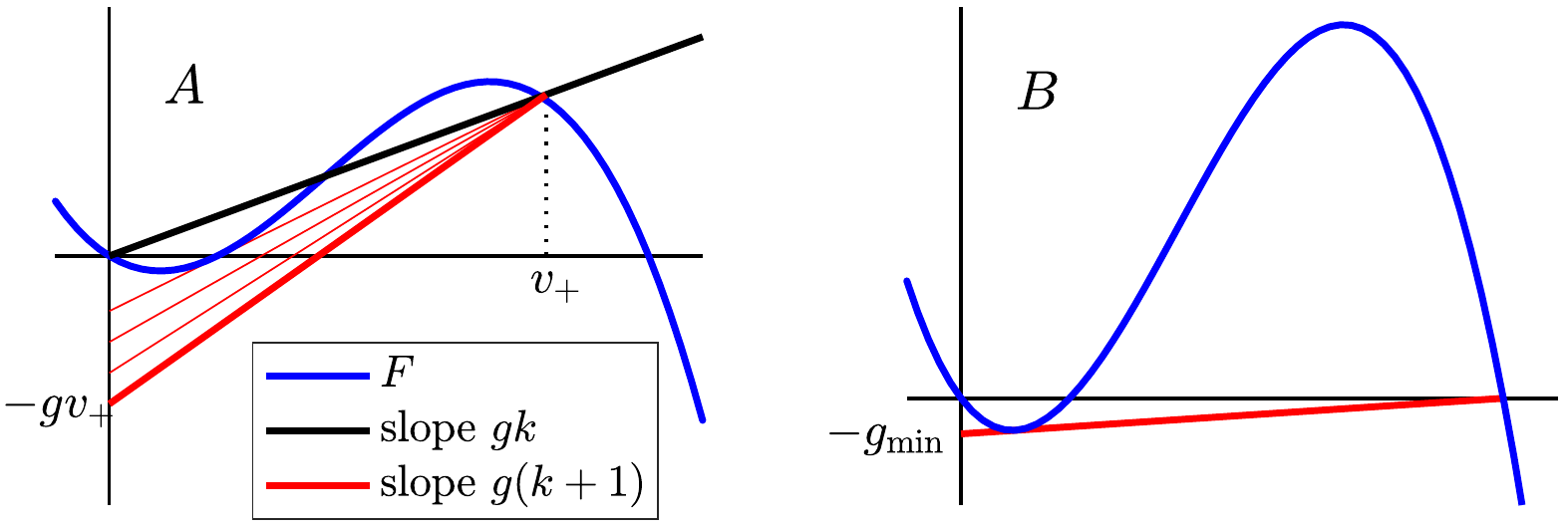}
\caption{A: Illustration of the condition in Theorem \ref{theorem:MK}. The red line must stay below the critical segment of the cubic or, at most, touch it tangentially for the limit of the $v_j$ to be $v_+$. Otherwise, iterates of $\varphi$ will eventually jump to the critical segment and then converge to 0. 
B: Since the vertical intercept of $G$ is $-gv_+$, the minimal value of $g$ for which there can be propagation is obtained by setting $k=0$, so that $v_+=1$. We then draw the tangent to the critical segment passing through $(1,0)$. This tangent intersects the vertical axis at $-g_{\rm min}$.}
\label{fig:MUNRO_KRULL}
\end{figure}

Because of the discontinuity, we won't use the standard theory of iterated one-dimensional maps, but explicitly, geometrically construct the iterates. We start with $v_0=1$, and define $v_1,v_2,\ldots$ by $v_{j+1} = \varphi(v_j)$, $j=0,1,2,\ldots$. Given $v_j$, $v_{j+1}$ is the smallest solution of $F(v) + g (v_j - v) - gk v = 0,$ that is, of 
$$
F(v) = g(k+1) v - g v_j.
$$
The right-hand side of this equation defines a straight line with slope $g(k+1)$ that passes through the point $(v_j, g k v_j)$. The construction of $v_{j+1}$ from $v_j$ is illustrated by Fig.\ \ref{fig:CONSTRUCTING_PHI}A: Given $v_j$, we draw the line with slope $g(k+1)$ through the point $(v_j, gk v_j)$. The $v$-coordinate of its smallest intersection with the cubic is $v_{j+1}$. From the figure, we read off the following result.  
\begin{theorem} 
\label{theorem:MK} Define $v_E$ as in Section \ref{sec:neighbors_at_rest} (shown in Fig. \ref{fig:CONSTRUCTING_PHI}B). 
If $gk \leq F'(v_E)$, we denote the positive solutions of the equation $F(v) = gkv$ by $v_- \in [v_T, v_E]$ and $v_+ \in [v_E,1]$; see Fig.\ \ref{fig:CONSTRUCTING_PHI}B. Starting with $v_0=1$, then $\lim_{j \rightarrow \infty} v_j$ is either $v_+$ or $0$. It is $v_+$ if and only if the line through $(v_+, F(v_+))$ with slope $g(k+1)$ lies below  the critical segment of $F(v)$, or at most touches it tangentially; see Fig.\ \ref{fig:MUNRO_KRULL}A. If $gk > F'(v_E)$, then $\lim_{j \rightarrow \infty} v_j = 0$. 
\end{theorem}

\subsection{The region of persistent propagation} \label{sec:per-prop} 

We say that {\em there is persistent propagation} if the sequence $\{v_j\}$ converges to $v_+$. The geometric construction described in Fig.\ \ref{fig:CONSTRUCTING_PHI} implies that there can only be persistent propagation if the graph of $G$ is below or tangent to the critical segment, as illustrated in Fig.\ \ref{fig:MUNRO_KRULL}A. Therefore, persistent propagation is possible if $g \geq g_{\rm min}$, where $g_{\rm min}$ is defined in Fig.\ \ref{fig:MUNRO_KRULL}B. Note, for $g = g_{\rm min}$, there is propagation only for $k = 0$. 
 
\begin{figure}[ht!]
\centering
\includegraphics[scale=0.78]{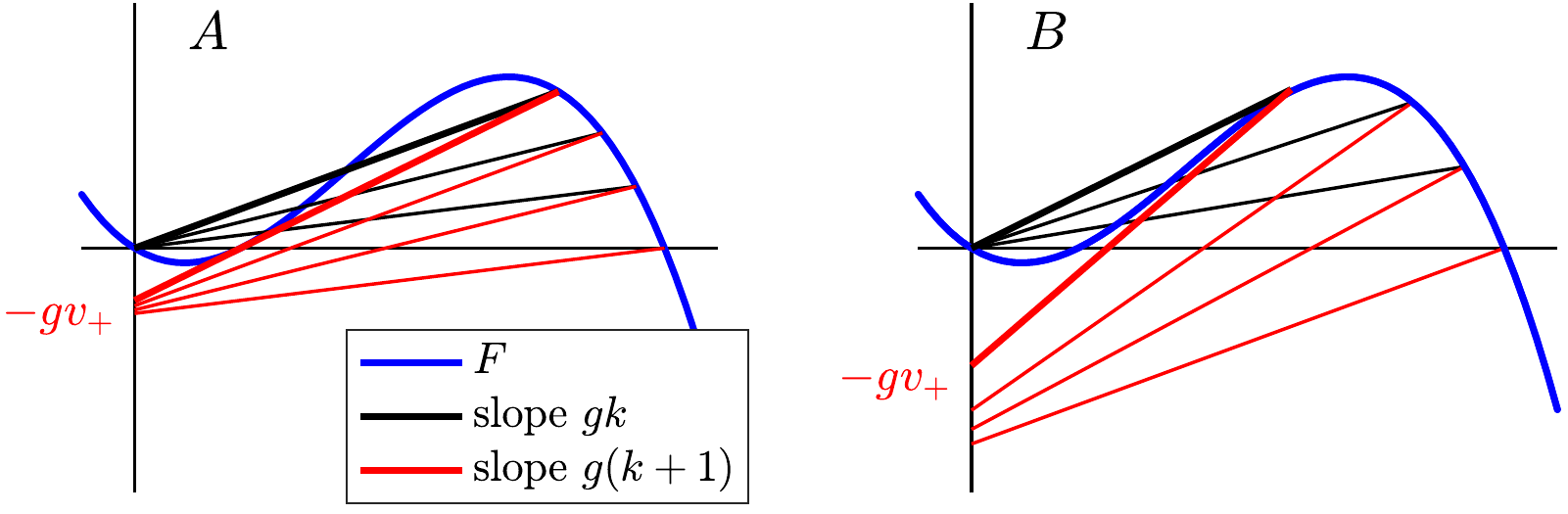}
\caption{A: For small $g$, the largest $k$ for which there is persistent propagation is determined by the condition that the line through $(v_+,gkv_+)$ with slope $g(k+1)$ touches the cubic tangentially (highlighted in bold). B: For large $g$, the largest $k$ for which there is persistent propagation is determined by the condition $gk = F'(v_E)$ (highlighted in bold).
}\label{fig:RANGE}
\end{figure}

For a given $g > g_{\rm min}$, we obtain the range of $k$-values for which there is propagation by watching how the lines $L(v)=gkv$ and $G(v)=g(k+1)v-gv_+$ rotate as $k$ increases, starting at $k=0$. Examples are shown in Fig.\ \ref{fig:RANGE}.
 
From Fig.\  \ref{fig:RANGE}, we read off the following result. 

\vskip 10pt
\begin{theorem} The region of persistent propagation in the $(g,k)$-plane is defined by
inequalities of the form $g \geq g_{\rm min},  0 \leq k \leq k_{\rm prop}(g)$, where $g_{\rm min}$ is defined by Fig.\ \ref{fig:MUNRO_KRULL}B, and $k_{\rm prop}$  is a continuous function with $k_{\rm prop}(g_{\rm min}) = 0$. For small $g$, $k_{\rm prop}$ is defined by the tangency condition illustrated in Fig.\ \ref{fig:MUNRO_KRULL}A. For large $g$, $k_{\rm prop}(g) = F'(v_E)/g$.
\end{theorem}
 
\begin{figure}[ht!]
\centering
\includegraphics[scale=0.5]{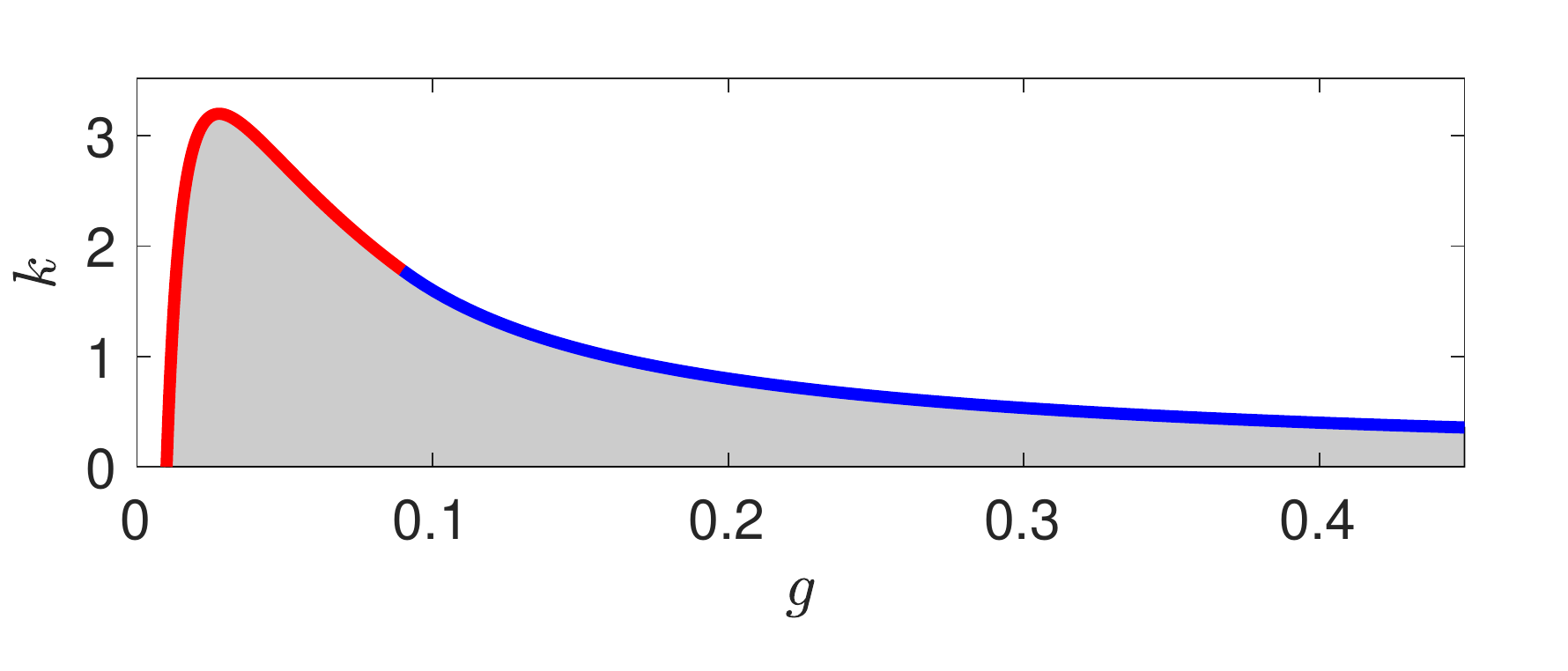}
\caption{The region of propagation in the $(g,k)$-plane. The red part of the curve indicates where $k_{\rm prop}(g)$ is determined by the tangency condition; the blue part indicates where $k_{\rm prop}(g) = F'(v_E)/g$.}
\label{fig:PROPAGATION_REGION}
\end{figure}

\vfill
\pagebreak
Since $k_{\rm prop}$ is $0$ at $g_{\rm min}$, and tends to zero as $g \rightarrow \infty$, there is  value of $g$ that is ``optimal" for propagation in the sense that $k_{\rm prop}$ is maximal. Figure \ref{fig:PROPAGATION_REGION} shows the region of propagation for one particular choice of parameters. The exact location of the peak depends on $v_T$, as it did in the single cell (see eq.\ (\ref{peak_equation})). 
 
\section{Comparison with more realistic models}\label{sec:realistic_model}

We apply our analysis to more realistic models of a neuron and  a myocyte. 

\subsection{A neuronal model}\label{sec:realistic_neuron}

We compare the single cell analysis of Section \ref{sec:single_cell} with simulation results for a Hodgkin-Huxley type model of the neuronal axon \cite{Traub:1999uf}. We consider a single cell with downstream neighbors held at rest (which is $v=0$ in \cite{Traub:1999uf}). To apply our analysis, we reduce the Hodgkin-Huxley type model to one dimension by replacing $m$ with $m_{\infty}(v)$, and keeping $n$ and $h$ fixed at resting values. Using the one-dimensional model we calculate $v_F$, the peak voltage of the disconnected cell, and $v_E$, the maximum possible threshold voltage (see Fig. \ref{fig:ZERO_VU_2}). 

\begin{figure}[htbp]
\begin{center}
\begin{minipage}{2.5in}
A
\centerline{single cell}
\resizebox{2.5in}{!}{\includegraphics{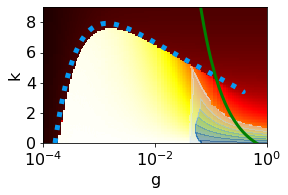}}
\end{minipage}
\begin{minipage}{2.5in}
B
\centerline{tree network}
\resizebox{2.4in}{!}{\includegraphics{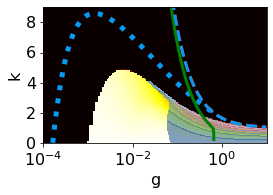}}
\end{minipage}

\begin{minipage}{2.5in}
C
\centerline{single cell}
\resizebox{2.5in}{!}{\includegraphics{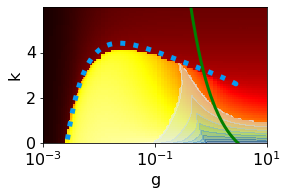}}
\end{minipage}
\begin{minipage}{2.5in}
D
\centerline{tree network}
\resizebox{2.4in}{!}{\includegraphics{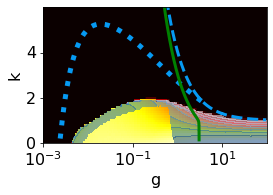}}
\end{minipage}
\caption{{\bf Behavioral regions for a Hodgkin-Huxley type neuronal model and the Luo-Rudy myocyte model.} (A) The heat map shows the maximum voltage for a single cell Hodgkin-Huxley neuron where the upstream cell's voltage is held at $v_u=v_F$. The change in maximum voltage, when we only hold $v_u=v_F$ until $v=v_E$, is shown with overlaid blue shading: The darker the shading the higher the difference. No shading indicates no difference between conditions. The dashed blue curve shows $k_{\rm max}$. The solid green curve shows the predicted active-passive boundary $k_{\rm exc}$. 
(B) The heat map shows the maximum voltage for the penultimate cell in the tree network where $v_0=v_F$ for the entire simulation. Overlaid blue shading shows the difference in maximum voltage where $v_0=v_F$ only until $v_1=v_E$. The curves $k_{\rm prop}$ (dashed blue) and $k_{\rm exc}$ (solid green) are both adjusted for the assumption that the downstream cell's voltage is proportional to the cell's voltage. 
(C+D) Single cell and tree network results for the Luo-Rudy myocyte model.
\label{fig:HHfull}}
\end{center}
\end{figure}

We run simulations for many $(g,k)$ pairs under two different conditions. In condition 1, the upstream cell is held at $v=v_F$ for the entire simulation. In condition 2, the upstream cell is set to $v=v_F$ until the central cell's voltage reaches $v_E$. For each simulation, we calculate the maximum voltage reached by the central cell, along with the difference in maximum voltages between conditions 1 and 2. If the central cell fires an AP under both conditions, then we say that active propagation is possible. (A precise definition of ``fires an AP" is not needed here because the results we present are quantitative, not qualitative.) However, if the central cell only fires an AP under the first condition, then only semi-active propagation is possible since it cannot fire a full AP on its own. These results are shown in Fig. \ref{fig:HHfull} panel A.

We then calculate $k_{\rm max}(g)$ and $k_{\rm exc}(g)$ using the one-dimensional reduction, and compare them directly to the simulation results. Figure \ref{fig:HHfull} panel A shows that $k_{\rm max}(g)$ for a single cell matches the active propagation region quite well. The boundary for semi-active propagation is not quantitatively accurate, but lies to the left of $k_{\rm exc}(g)$. The quantitative difference between the predicted and actual results is not surprising, considering that our predictions are based on a highly simplified model in which, for instance, n and h stay fixed at resting values.

For the collapsed tree network, we make modifications in both the analysis and the simulations. In the simulations, we only hold downstream neighbors fixed at $v=0$ for the last cell. In our analysis, we therefore adjust the assumption that downstream neighbors are held at $v=0$. Instead, we make a cell's downstream neighbor voltages proportional to the cell's voltage \cite{Ramasamy:2007bz}, which follows the steady-state condition when there are only passive currents. The details of this modification are explained in the appendix.

In Fig.\ \ref{fig:HHfull} panel B, we see that the active propagation region for the tree network lies entirely within the predicted region and has a similar shape. That $k_{prop}$ over-estimates the boundary for AP propagation is to be expected, since our analysis does not take the duration of an AP into account. As before, the semi-active propagation region also extends to the left of $k_{\rm exc}(g)$. 

\subsection{A myocyte model}

We reduce the Luo-Rudy myocyte model \cite{Luo:1991vz} to one dimension by replacing $m$ with $m_{\infty}(v)$ and holding all other gating variables fixed at their resting values. We chose how to replace the gating variables by comparing their time constants with the membrane time constant when the cell is at rest. For ease of calculation, we also re-center $F(v)$ so that the resting fixed point is $v=0$. Using the same simulation setup as in section \ref{sec:realistic_neuron}, we again see that propagation into a single cell matches the predicted boundaries quite well. Likewise, the predicted active propagation region for the tree surrounds the active propagation region in the simulations. The semi-active propagation region boundary lies to the left of the predicted boundary.

\overfullrule=0pt

\section{Discussion}

We explained why AP propagation through a network of excitable cells connected by gap junctions depends non-monotonically on gap junction strength. We did this by reducing models to one dimensional firing currents, which allows us to visualize how the firing currents and gap junction currents  interact directly. We also found a new behavioral region where cells may propagate APs but are so strongly connected that they are no longer individually excitable. 

We believe this framework can give a simple and efficient method for understanding where AP propagation may occur. Moreover, this framework may be used to understand how network behavior may change with changes in connectivity, gap junction strength, and firing currents. While we did not model a gap junction network with heterogenous connections, our framework may still apply, since often AP propagation is determined by a make-or-break point with the highest connectivity. We may also predict when heterogenous connectivity may lead to re-entrant or spiral waves within a gap junctional network \cite{Munro:2010ih}. Finally, while we focused on AP propagation here, our framework also offers a way of visualizing subthreshold network characteristics, including how gap junction strength affects the firing threshold and how activation currents mix with gap junction currents to determine a subthreshold response. These characteristics may play a role in network synchronization mediated by gap junctions.

AP propagation through gap junctions is seen in many contexts. This means our framework may apply to a variety of cell types and networks. For example, gap junctions mediate heart contractions \cite{Bernstein:2006, Carmeliet:2019, Kieval:1992wg}, and possibly other muscle contractions \cite{Borysova:2018}. APs actively propagate across neuronal gap junctions throughout the nervous system, including networks of interneurons and axons \cite{Shimizu:2013}. Gap junctions also mediate propagating calcium waves \cite{Goldberg:2010, Kahne:2019, Malmersjo:2014}. While the cells involved may not exhibit full APs, the cells' partial excitability may be analyzed in a similar manner. Lastly, gap junction conductance is highly plastic, and changes in conductance may result in change of function \cite{Haas:2011fk, OBrien:2018, Pereda:2004cl, Szabo:2010bk}, or maintain a similar function as connected cells change over time \cite{Parker:2009cz}.

There are also many medical conditions where either gap junction conductance or firing currents of gap junctionally connected cells are affected. Cardiac arrhythmia can be brought about by pathogenic stressors \cite{Calhoun:2020, Kim:2013}, adrenergic stimulation \cite{Salameh:2011jm}, or myocardial ischemia (buildup of plaques) associated with decreased gap junctional coupling \cite{Dhein:2006ej, Kieval:1992wg}. Abnormal gap junction ubiquitination can affect both the heart and nervous system \cite{Totland:2020}. Gap junctions may be sensitive to toxins \cite{Vitale:2018}. Epilepsy has been tied to gap junctions connecting networks of pyramidal cell axons, interneurons, and astrocytes, each of which play a different role in network excitability \cite{Carlen:2000vt, Lapato:2018, Szente:2002ul, Traub:2011}. Furthermore, gap junctions are indicated in Parkinson's disease \cite{Schwab:2014}, Alzheimer's disease \cite{Wang:2009fw}, nerve injury \cite{Chang:2000tw, Murphy:1983uc, Nagaoka:1999ux}, and degenerative diseases in the retina \cite{OBrien:2018}. Modulating gap junction strength may help with several of these conditions, however questions remain on how much gap junctions can be modulated without producing adverse side effects \cite{Li:2019, Manjarrez-Marmolojo:2016}. An alternative may be to modulate firing currents through gene therapy \cite{Hucker:2017}. There are also many conditions where external stimulation is applied as a part of therapeutic treatment \cite{Bikson:2004km, Dell:2019, Prakosa:2018, Trayanova:2014, Zangiabadi:2019}. We hope that our framework can give some insight on appropriate ranges to maintain proper function in a variety of networks. 

There are numerous studies using model networks and simulations to tease out the mechanisms and characteristics of different gap junctional networks. Predicting conduction through the heart is an active line of research, where experiments and simulations work hand in hand \cite{Arutunyan:2003cp, Clayton:2005er, Hand:2010, Henriquez:2001ia, Hubbard:2007ka, Hurtado:2020, Keener:1998ke, King:2021, MacCannell:2007eo, Nygren:2000, Priya:2015, Qu:2004kx, Ramasamy:2007bz, Spach:1995uta, Spach:2000wqa, Steinberg:2006em}. Researchers also use simulations to study active propagation across neuronal dendrites \cite{Gansert:2007gm, Nadim:2006eq} and axons \cite{Cunningham:2012co, Lewis:2000uf, Lewis:2003uo, Maex:2007bi, Munro:2010ih, Munro:2012iy, SteynRoss:2007ic, Traub:1999uf}. Work on propagating calcium waves, where signals can actively propagate both intracellularly and between cells by gap junctions, uses a mixture of theory and model simulations \cite{Dougoud:2016, Harris:2020, Keizer:1998}. Active propagation seen in gap junction networks also relates to active propagation along neural axons, where model simulations have different cell compartments connected linearly following Ohm's law \cite{Alcami:2019, Hull:2015, Tuckwell1988}.

We note that our framework only provides simple guidelines for outlining where AP propagation may occur. While this theory can give a first estimate on when to expect propagation, there are many ways we could increase the accuracy of this prediction. For instance, we could link different versions of $\varphi(v)$ to model propagation through a network with heterogenous numbers of downstream neighbors. We may also improve the prediction by taking the connectivity of downstream neighbors more fully into account. In fact, propagation may not only depend on the number of immediate downstream neighbors, but the number of neighbors several steps away may affect propagation \cite{Munro:2010ih}. The predicted AP propagation regions in our framework overestimate the regions found in the original models. This may be due to higher-dimensional interactions within the cell as well as how long an AP lasts vs.\ the time it takes an AP to trigger an AP in a downstream cell. Extending our analysis to include the interaction of recovery currents with the firing currents may improve our prediction. Likewise, taking timing into account by studying AP duration vs.\ the separation between the graphs of $F$ and $G$ may also improve prediction, especially of AP propagation through a network. 

Our framework studies propagation assuming the first cell is firing. It doesn't address getting the first cell to fire in the first place, but rather reproduces the situation where the cell is voltage clamped. As gap junction conductance increases, it may take an increasing amount of current to bring the first cell to the desired voltage. In some cases, we may need to depolarize multiple downstream neighbors in the network to get AP propagation, similar to what is seen in semi-active propagation.

Overall, reduction of cell models to one-dimensional firing currents gives insight into how gap junction currents can allow active propagation, semi-active propagation, or passive propagation. The shape of the firing currents explains how the threshold and firing voltage change with conductance, and allows us to predict an ideal gap junction conductance for propagation. This framework gives further insight into the mechanisms of these networks - which play a key role in both normal biophysiological function and disease. 

\section*{Acknowledgments} We thank Nancy Kopell for many insightful comments on this work.
~
\vskip 30pt

\bibliographystyle{siam} 
\bibliography{refs.bib, gj_prop.bib}

\vskip 30pt
\noindent
{\bf Appendix. Relaxing the assumption of downstream neighbors at rest.}
\vskip 5pt
Throughout our analysis, we assumed that each cell responds to its upstream neighbor, but behaves as if its downstream neighbor stayed at rest, $v=0$. However, we modified this assumption in section \ref{sec:realistic_model} so that the downstream neighbor's voltage is proportional to the cell's voltage as suggested in \cite{Ramasamy:2007bz}. 

To explain the modified assumption, think about an infinite chain of cells with passive current and voltages $v_j$, $j=0,1,2,\ldots$. The $j$-th cell ($j \geq 1$) sees the upstream voltage $v_{j-1}$, and the downstream voltage $v_{j+1}$. If the voltages are in equilibrium, then 
\begin{equation}
\label{eq:steady_state}
g (v_{j-1} - v_j) + gk (v_{j+1} - v_j) - g_L v_j = 0, ~~~~ j = 1, 2, 3, \ldots.
\end{equation}
This is a linear difference equation. From the standard theory of such equations, the solutions are the sequences
\begin{equation}
\label{eq:chain_solution}
v_j = C_+ \alpha_+^j + C_- \alpha_-^j,
\end{equation}
where $C_+$ and $C_-$ are constants, and $\alpha_+$, $\alpha_-$ are the solutions of the quadratic equation 
$$
g \left( 1 - \alpha \right)  + gk (\alpha^2 - \alpha) - g_L \alpha = 0, 
$$
or equivalently, 
\begin{equation}
\label{eq:quadr}
 k\alpha^2 -  \left( k + 1 + \frac{g_L}{g} \right)   \alpha + 1= 0.
\end{equation}
It is convenient to write 
$$
\beta = 1 + \frac{g_L}{g}.
$$
The solutions of (\ref{eq:quadr}) are then
\begin{equation}
\label{alphapm}
\alpha_{\pm}  =  \frac{k + \beta \pm \sqrt{ \left( k + \beta\right)^2 - 4k}}{2k}.
\end{equation}
It is straightforward to verify that
$$
0 < \alpha_- < 1 < \alpha_+.
$$
Therefore the only constant {\em bounded} solutions of (\ref{eq:steady_state}) are 
$$
C \alpha^j
$$
where $C$ is constant and
$$
\alpha=  \alpha_- =  \frac{k + \beta - \sqrt{ \left( k + \beta\right)^2 - 4k}}{2k}.
$$
The cell downstream to the $j$-th cell is at voltage $v_{j+1} = \alpha v_j$.

The modified assumption underlying the analysis results of section \ref{sec:realistic_model} is that a cell at voltage $v$ sees the downstream voltage 
$$v_d = \alpha v$$ 
{\em always}, not just at equilibrium. This modification changes our earlier analysis in a surprisingly simple way. The gap-junctional current 
$$
g \left( v_u - v \right)  - gk v
$$
is replaced by 
$$
g \left( v_u - v \right) + g k (v_d - v) = 
g \left( v_u - v \right)  + g k (\alpha v - v)  = g \left( v_u - v \right)  - g k (1 - \alpha) v.
$$
So what used to be ``$k$" is now ``$k (1-\alpha)$". For instance, the condition $k<k_0(g)$, where $k_0(g)$ is a boundary such as $k_{\rm exc}(g)$ from section \ref{sec:single_cell}, becomes 
$$
k \left( 1 - \alpha \right) < k_0(g).
$$
This is equivalent to 
\begin{equation}
\label{eq:new_ineq}
\frac{k - \beta +\sqrt{  \left( k + \beta\right)^2 - 4k } }{2} < k_0(g).
\end{equation}
Using the definition of $\beta$, we find that (\ref{eq:new_ineq}) is equivalent to 
$$
k <   k_0(g) \left(1 + \frac{1}{k_0(g) +  \frac{g_L}{g}}\right).
$$

\end{document}